# On Generic Frequency Decomposition.[1]
# Part 1: Vectorial decomposition

**Sossio Vergara**
**Universidad ORT Montevideo, Uruguay**
**ITI B. Pascal Rome, Italy**
sossio@montevideo.com.uy

## *Abstract*

The famous Fourier theorem states that, under some restrictions, any periodic function (or real world signal) can be obtained as a sum of sinusoids, and hence, a technique exists for decomposing a signal into its sinusoidal components. From this theory an entire branch of research has flourished: from the Short-Time or Windowed Fourier Transform to the Wavelets, the Frames, and lately the Generic Frequency Analysis.
The aim of this paper is to take the Frequency Analysis a step further.
It will be shown that keeping the same reconstruction algorithm as the Fourier Theorem but changing to a new computing method for the analysis phase allows the generalization of the Fourier Theorem to a large class of nonorthogonal bases.
New methods and algorithms can be employed in function decomposition on such generic bases. It will be shown that these algorithms are a generalization of the Fourier analysis, i.e. they are reduced to the familiar Fourier tools when using orthogonal bases. The differences between this tool and the wavelets and frames theories will be discussed. Examples of analysis and reconstruction of functions using the given algorithms and nonorthogonal bases will be given.
In this first part the focus will be on vectorial decomposition, while the second part will be on phased decomposition. The phased decomposition thanks to a single function basis has many interesting consequences and applications.

Keywords: Fourier theorem, functional analysis, wavelets, frames, nonorthogonal basis.

---

## Introduction

J. B. J. Fourier enunciated his celebrated theorem at the beginning of the XIX century while he was looking for a way to solve the equation of the heat transfer in solids. The idea was that, as the differential equations governing the process were too difficult to solve, one could consider the result as the superposition of the effects of simpler functions, for instance the sinusoids; a technique still in use.

His original paper, dated 1807, was initially rejected due to the opposition of Lagrange, a great Italian mathematician, who thought that Fourier's claim that "any continuous and periodic function can be obtained as a sum of properly chosen sinusoidal waves" could not be applicable to signal with discontinuities. Finally the paper was published, but not before Lagrange's death, many years later.

It now turns out that both Fourier and Lagrange were partially correct, but some advances in the field of mathematics were necessary in order to give the Theorem its exact enunciation.

As is well known, the Theorem in its standard form states that any function (and hence any real world signal) that is periodic, has finite energy and a limited number of discontinuities in a cycle, can be decomposed as a series of sine and cosine or other orthogonal functions.

Such functions are called the bases of that space, and the series is said to converge to the function. Actually, some different definitions of convergence exist: pointwise, uniform and $L^2$-norm, the last being the most satisfactory from the mathematical point of view.

By definition the norm of a function $g(x)$, periodic with a period P, is:

$$\| g(x) \| = \left[ \frac{1}{P} \int_0^P | g(x) |^2 \, dx \right]^{\frac{1}{2}} \tag{1}$$

And the series of partial sums $S_N$ is said to converge to $g(x)$ in norm if:

$$\lim_{N \to \infty} \frac{1}{P} \int_0^P | g(x) - S_N(x) |^2 \, dx = 0 \tag{2}$$

More details can be found in [6], [9] and [10].

In the last two centuries, this theorem has fueled research in the field of pure mathematics as well as in physics. From the original series, an entire set of tools has been developed, including the Transform, then adapted to the DFT, that finally evolved into its speedy version, the FFT.

Given its importance, much work has been dedicated to the quest for a generalization of the Fourier Theorem, primarily in the direction of searching for other orthogonal bases, like the Legendre polynomials [6], [10].

A limitation of the standard Fourier transform is its inability to discriminate signals in time, as the analysis is always performed on the entire duration of a signal. For example, the information on any abrupt change in the signal is spread out on the entire frequency axis.

To overcome this limitation another mathematical tool has been developed: the Short Time Fourier Transform (STFT). In the last decades a new tool has evolved: the wavelets [1][2][8]. These, initially developed as orthogonal functions [5], were then generalized to the non-orthogonal case. They have the advantage over the Fourier transform of being localized both in time and frequency, while changing the window duration with the frequency: it is a multiresolution analysis.

In wavelet terminology the concept of scale is preferred to frequency. It is a time-scale representation of the signals.

Unfortunately, the wavelets have some disadvantages in the reconstruction phase [7]; moreover, only a limited set of wavelets is available and they all are special mathematical functions [1], [2].

In the last few years a new concept has emerged as an abstract framework for the time-frequency decomposition: the Frames [4] [3] [1]; this theory analyzes the completeness, stability and redundancy of discrete signal representations. The frame theory establishes the general conditions under which one can recover a vector *f* in a Hilbert space **H** from its inner products with a family of vectors $\{f_n\}_{n\in\Gamma}$ called the "frame" ( $\Gamma$ being an index set ).

However the frame theory still has some limitations: to reconstruct a signal from its inner products one needs to invert the frame operator and this is not an easy task. In some cases this can be carried out only employing a dual frame. Due to these difficulties, only few special types of frames are actually employed, and mainly in the transmission and reconstruction of signals in noisy environments.

Lately, a different approach to the generalization of the Fourier Theorem has been proposed: the Generic Frequency Analysis, [11], [12], [13]. It has been demonstrated that any function *f* in $L^2$ can be reconstructed by using even and odd couples of functions that have multiplicative Fourier coefficients; while the problem of finding the general rules of the frequency analysis is considered still unsolved [14]. These results are quite recent and applications for them have not yet been found. Indeed they limit the analysis to the use of even and odd couples of specific functions like the square waves, the trapezoidal and the like.

There are evident differences between the Generic Frequency Analysis and the wavelets and frames, as these latter have been developed especially for localized time-scale analysis. A generalization of the Fourier Theorem based on generic functions could have some advantages and many important implications.

These advantages will be evident if we can loosen the requisites for a couple of functions in order to be a basis in $L^2$; and even more if we find a way to decompose a function in terms of such a generic basis but using polar coordinates, possibly avoiding the use of dual bases for analysis and reconstruction.

The aim of this paper is to show that such generalization of the frequency analysis is possible, provide the necessary tools and to indicate some of the applications and consequences.

We first need to prove that any function or couple of functions that satisfies some quite loose requisites can be a basis in the space of the functions for which the Fourier Theorem is valid (the $L^2$ Hilbert space). In other words: given two functions $S(x), R(x) \in L^2$, satisfying those requisites, for any function $f(x) \in L^2$, exists a unique set of coefficients $A_k$ and $B_k$, such that the series of $A_k S(kx) + B_k R(kx)$ converges to *f(x)* in norm. So it can be written:

$$f(x) = C_0 + \sum_{k=1}^{\infty}[A_k S(kx) + B_k R(kx)] \qquad (3)$$

Similarly, when using polar coordinates any function $f(x) \in L^2$ can be reconstructed univocally as a superposition of $S(kx)$ as:

$$f(x) = C_0 + \sum_{k=1}^{\infty} M_k S(kx + \Phi_k) \qquad (4)$$

And hence the function $S(kx)$ or the couple $[S(kx), R(kx)]$ is a basis in $L^2$.

Some computationally simple algorithms will be developed for the analysis. It will be shown that these algorithms turn out to be the standard Fourier tool when orthogonal bases are employed.

The application of the new theory in building a "generalized" form of DFT is straightforward. Some practical applications have already been developed, and other, that can promptly benefit from the theory, are under investigation.

In this first part of the paper we shall examine the case of vectorial decomposition, while in the next part we shall study the case of the polar coordinates together with its consequences and applications.

Throughout the paper the emphasis will be on intuitive comprehension and applications, as the goal of the paper is to explore to the greatest extent the possibilities of the new tool, even at the cost of sacrificing some details for the moment. As such when it is straightforward, generalization will be omitted in favor of a more direct approach to the results; for example, real valued functions and series of functions will be treated for the most part, as the case of complex functions and transforms can be easily derived. And, when not otherwise specified, the period P of the functions is [0, 1] so to neglect this constant in the equations.

## *The classic parallelism between vector and function spaces*

It is well known that a function *f(x)*, periodic with period P = 2π, that satisfies the Dirichlet conditions:
1) *f(x)* has a finite number of maxima and minima in one period.
2) *f(x)* has a finite number of discontinuities in one period.
3) and: $\int_0^P |f(x)| dx < \infty$

can be decomposed in a series of sine and cosine:

$$f(x) = C_0 + \sum_{k=1}^{\infty} (a_k \sin kx + b_k \cos kx) \qquad (5)$$

Where:

$$C_0 = \frac{1}{2\pi} \int_{-\pi}^{+\pi} f(x) dx \qquad (6)$$

$$b_k = \frac{1}{\pi} \int_{-\pi}^{+\pi} f(x) \cos kx \, dx \qquad (7)$$

$$a_k = \frac{1}{\pi} \int_{-\pi}^{+\pi} f(x) \sin kx \, dx \qquad (8)$$

Or, in the more compact, complex exponential form:

$$f(x) = \sum_{k=-\infty}^{\infty} c_k e^{jkx} \qquad (9)$$

Where:

$$c_k = \frac{1}{2\pi} \int_{-\pi}^{+\pi} f(x) e^{-jkx} dx \qquad (10)$$

The above is called the Fourier Series; it can be extrapolated to continuous and aperiodic functions by means of the Fourier Transform (and Discrete Fourier Transform in the case of Discrete-Time signals) and its inverse:

$$T_F(f(x)) = F(\omega) = \int_{-\infty}^{+\infty} f(x) e^{-j\omega x} dx \qquad (11)$$

$$f(x) = T_F^{-1}(F(\omega)) = \frac{1}{2\pi} \int_{-\infty}^{+\infty} F(\omega) e^{j\omega x} dx \qquad (12)$$

The sufficient condition for the existence of the Fourier Transform is that *f(x)* be (Lebesgue) square integrable:

$$\int_{-\infty}^{+\infty} |f(x)|^2 dx < \infty \qquad (13)$$

or equivalently *f(x)* ∈ L$^2$.

When studying these series, the parallel between the function and vector spaces is often made, showing that a function, as a vector, can be decomposed in its projections on two orthogonal axes. In the vector space these axes are **x** and **y**, in the function space they are the sine and cosine, or the real and imaginary part.
But let us see in some more detail what happens in the vector space.
Can we decompose a vector in its components along a new system of reference whose axes are not perpendicular to each other? Is this decomposition unique?
Let us make a drawing in the plane as in Fig. 1.

Here x, y is our standard orthogonal Cartesian system of reference; S and R is a new non-orthogonal system of reference, and V is a generic vector. It is easy to see that any vector V can be univocally decomposed in the non-orthogonal system of reference formed by S and R axes.
But only on condition that the ratios of the modules of projections of the unity vector of the new R and S axes along x and y be different:
$|S_x| / |S_y| \neq |R_x| / |R_y|$
If these ratios are equal it means that R and S lay on the same line, so no decomposition is possible (in mathematical terms these two axes need to be independent, but not necessarily orthogonal, to span the entire space).
Now, can we extrapolate these findings to the function space? Could it be useful?
A first level extrapolation is quite trivial.
Having two functions as:

$$S(x) = s_1 \cos(x) + s_1' \sin(x) \qquad (14)$$

$$R(x) = r_1 \cos(x) + r_1' \sin(x) \qquad (15)$$

here the coefficients $s_1$, $s_1'$, $r_1$, $r_1'$ are the "projections" of $S(x)$ and $R(x)$ over the sine and cosine axes. If they are independent:

$$\left|\frac{s_1}{s_1'}\right| \neq \left|\frac{r_1}{r_1'}\right| \qquad \text{or:} \qquad |s_1 \, r_1'| \neq |s_1' \, r_1| \qquad (16)$$

any function that can be decomposed in a Fourier series, can also be written as a series of the functions $S(x)$ and $R(x)$ as:

$$f(x) = C_0 + \sum_{k=1}^{\infty} [A_k S(kx) + B_k R(kx)] \qquad (17)$$

The proof is trivial and will be omitted.
It should be noted that these new basis functions are non-orthogonal to each other.
Two periodic functions $f_1(x)$ and $f_2(x)$ are said to be orthogonal if they have a null inner product:

$$\langle f_1(x), f_2(x) \rangle = \int_0^P f_1(x) \overline{f_2(x)} dx = 0 \qquad (18)$$

Here the second function is complex conjugated, P being the period.
In the case of the previously defined functions $S(x)$ and $R(x)$ we have:
$\langle S(kx), R(kx) \rangle \neq 0 \quad \forall k$ (19)
but:
$\langle S(kx), S(ix) \rangle = 0 \quad \forall k \neq i$ (20)

We now introduce two definitions:
we shall call "horizontally orthogonal" a couple of periodic functions $f_1(x)$, $f_2(x)$ such that :
$\langle f_1(kx), f_2(kx) \rangle = 0 \quad \forall k$

We shall call "vertically orthogonal" a couple of functions if:

$<f_1(kx),f_1(mx)> = 0$ $\quad\quad$ $<f_1(kx),f_2(mx)> = 0$ $\quad\quad$ $<f_2(kx),f_2(mx)> = 0$ $\quad \forall\ k \neq m$

Then the functions S(x) (14) and R(x) (15) as defined above, form a couple of functions that are horizontally non-orthogonal but vertically orthogonal.
Unfortunately not much can be done with such simple functions that are, essentially, sinusoids out of phase.
So let us see if we can extrapolate these results to some, more generic, class of functions.
A further generalization has been carried out in [11], [12]. There it was demonstrated that couples of even and odd functions that have multiplicative Fourier coefficients are a basis in $L^2$. These functions are now horizontally orthogonal but vertically non-orthogonal. Details can be found in [12].
This analysis procedure starts by choosing a suitable couple of functions X(x) and Y(x); from these, the family of their biorthogonal functions: $h_n(x)$, $g_n(x)$ is generated. Finally, the coefficients are calculated via the inner product of the function $f(x)$ with any of these $h_n(x)$, $g_n(x)$.
The main restriction here (and in the previous tools as the frames) seems to be the necessity of using the inner product in the computing of the coefficients; this, in turn, forces the use of biorthogonal functions and limits the set of functions that can be used as bases to even and odd couples of square waves and few others.
To overcome these limitations we now propose a new paradigm for the function decomposition, such as: *"The analysis is anything that is needed to find the coefficients for the reconstruction of a function."*
Well, it seems a tautology or at least trivial (as, at first look, some propositions of philosophy [15]); but the change here is that reconstruction and analysis are now undefined. They are just linked to each other.
As a matter of fact, different reconstruction algorithms could be devised, all of them capable of approximating any function $f(x)$ starting from some basis. Any of these reconstruction algorithms needing its own specially computed set of coefficients.
Of course one could better introduce the concept of limit, the requirement of uniqueness of the reconstruction and hence the basis, but let us keep the statement as simple as it is for the moment.

Now, if the postulate is acceptable, we can develop a new analysis/reconstruction structure, whose only constraint is the capability of univocally approximating any function. This structure shall not be limited to the use of the inner product in the computation of the coefficients. All we require is that increasing the order of approximation, the error (in norm) between the function itself and its reconstruction tends to zero, and that these coefficients are unique, so that we have a bijective relation between the space of functions and the space of coefficients, that in turn allows us to define an inverse to that relation.

It can be demonstrated that two functions:

$$S(x) = \sum_{i=1}^{\infty} \left[ s_i \cos(ix) + s_i' \sin(ix) \right] \quad\quad (21)$$

$$R(x) = \sum_{i=1}^{\infty} \left[ r_i \cos(ix) + r_i' \sin(ix) \right] \quad\quad (22)$$

such that:
$\quad\quad |s_1 r'_1| \neq |s'_1 r_1| \quad\quad (23)$

if for any $A_1$, $B_1$ not both zero

$$(A_1 s_1 + B_1 r_1)^2 + (A_1 s_1' + B_1 r_1')^2 > \sum_{i=2}^{\infty} (-A_1 s_i - B_1 r_i)^2 + \sum_{i=2}^{\infty} (-A_1 s_i' - B_1 r_i')^2 \quad\quad (24)$$

*Then:* for any $f(x) \in L^2$ exists a unique set of coefficients $C_0$, $A_k$, $B_k$, such that the series:

$$S_N = C_0 + \sum_{k=1}^{N} [A_k S(kx) + B_k R(kx)] \tag{25}$$

converges to *f(x)* in norm as $N \to \infty$.

A proof can be found in the Appendix.
Similarly, one can extend the previous result to polar coordinates decomposition. The issue will be examined in the second part of the article.
It can be verified that this tool is a superset of the generic frequency analysis as introduced in [12]; as a matter of fact, the requirements on the basis (equations 23 and 24 here) are quite more relaxed than those defined there.
Note how we base our work on the reconstruction algorithm only. The analysis will be functional to reconstruction.
Therefore, the previous assertion makes sense only if, given a basis and a reconstruction algorithm, we can find a method to calculate the coefficients $A_k$ and $B_k$.
So, let us develop two methods instead, one will be a direct, or brute force, approach to the computation of the coefficients, the other one will be indirect, in the sense that makes explicit use of a change of coordinates. This indirect method has some advantages over the previous one, as a simpler and faster iterative algorithm and the possibility of being adapted to polar coordinate decomposition.

### *"Direct" method*

The Eq. (25) for N components can be rewritten as (we consider for simplicity $C_0=0$):

$$S_N(x) = A_1 S(x) + A_2 S(2x) + \ldots + A_N S(Nx) + B_1 R(x) + B_2 R(2x) + \ldots + B_N R(Nx) \tag{26}$$

This will be the N-th approximation of the function *f(x)*. As $N \to \infty$, $S_N$ tends to *f(x)*, so we can substitute $S_N(x)$ with *f(x)* with an error decreasing as N increases. We can then write:

$$f(x) \cong A_1 S(x) + A_2 S(2x) + \ldots + A_N S(Nx) + B_1 R(x) + B_2 R(2x) + \ldots + B_N R(Nx) \tag{27}$$

We can now take the inner product of this equation with each one of the 2N functions:
$S(x), S(2x), \ldots S(Nx), R(x), R(2x), \ldots R(Nx)$.
Obtaining 2N equations.
Therefore, to find the 2N coefficients: $A_1, \ldots A_N, B_1, \ldots, B_N$, is sufficient to solve the system (28) of 2N equations:

$$(28)\begin{cases} A_1\int S(x)S(x)+...+A_N\int S(Nx)S(x)+B_1\int R(x)S(x)+...+B_N\int R(Nx)S(x)=\int f(x)S(x) \\ . \\ . \\ . \\ A_1\int S(x)S(Nx)+...+A_N\int S(Nx)S(Nx)+B_1\int R(x)S(Nx)+...+B_N\int R(Nx)S(Nx)=\int f(x)S(Nx) \\ A_1\int S(x)R(x)+...+A_N\int S(Nx)R(x)+B_1\int R(x)R(x)+...+B_N\int R(Nx)R(x)=\int f(x)R(x) \\ . \\ . \\ . \\ A_1\int S(x)R(Nx)+...+A_N\int S(Nx)R(Nx)+B_1\int R(x)R(Nx)+...+B_N\int R(Nx)R(Nx)=\int f(x)R(Nx) \end{cases}$$

There is no complex conjugation because we are dealing here with real valued functions or signals. The integrals are over a period.
Note that here the: S(kx)S(mx), R(kx)R(mx), S(kx)R(mx) integrals shall be comprised in the system only if they refer to a frequency that is included in the approximation; i.e.: only if km ≤ N, or if k is a factor of m.
The reason for this reduction of the system is obvious when one thinks of the functions S(kx) and R(mx) as Fourier series:

S(kx)=s'$_1$sin(kx) +s$_1$cos(kx)+ s'$_2$sin(2kx)+s$_2$cos(2kx)+s'$_3$sin(3kx)+…        (29)
R(mx)=r'$_1$sin(mx)+r$_1$cos(mx)+r'$_2$sin(2mx)+r$_2$cos(2mx)+r$_3$'sin(3mx)+…        (30)

The product of the functions S(kx)R(mx) (if the integral is nonzero) can have components only at frequencies that are present in both functions. The first common component of S(kx) and R(mx) will be at frequency mk or, eventually, at m but only if m is a multiple of k.
If those terms are not eliminated from the matrix (28), an error (and hence a noise) will be introduced due to the aliasing.
This reduction of the system makes the corresponding matrix quite sparse which helps in finding the solution. A further reduction of the matrix comes from the integrals that give zero as a result of vertical or horizontal orthogonality of the basis.
If a completely orthogonal base is employed, like the familiar sine-cosine, the system becomes purely diagonal, transforming itself into the renowned Fourier tool.

A software routine has been written to test the algorithm with different functions. As one could probably expect, the coefficients {$A_1$,..$A_N$, $B_1$,..$B_N$ } computed using this method are not exact and have the unpleasant behavior of changing with the number N of components. If we compute $A_1$,$A_2$, $B_1$, $B_2$ using in the system (28) five components, we get different results from the same $A_1$, $A_2$, $B_1$, $B_2$ when computed for six components.
The reason for this behavior is in the substitution we operated above, between $S_N(x)$ in Eq.: (26) and *f(x)* in Eq. (27) and in the system (28). As a result of the substitution, the system tries to find the lowest possible error between the reconstructed function and the *f(x)*, independently of the number of components N. As a matter of fact, when using a finite number of components, employing non–ideal coefficients can lead to a lower RMS error in reconstructing a function. Of course, increasing the number N of components, the coefficients calculated with this method tend to the theoretical ones.
Conversely, the subsequent indirect approach gives the exact reconstruction at each frequency thus generating the ideal coefficients.
Moreover, when it comes to the apparently overwhelming number of integrals needed, it should be noticed that the algorithm requires about $(2N)^2 + 2N$ integrals. But $(2N)^2$ of these, relative only to the basis, can be calculated "a priori"; a further reduction is possible by taking into account that, as we are here dealing with real valued functions the matrix is Hermitian:
< S(nx), S(kx) > = < S(kx), S(nx) >
The similarity between this approach and the frames theory is evident: the solution of the system (28) is equivalent to the inversion of the frame operator. Since the frames are such a comprehensive tool, any vectorial decomposition can be studied as a frame. Nevertheless, frames suffer from a general difficulty in finding the inverse operator that greatly limits their applicability, while the system (28) here can be seen as the product of a change of coordinates operator over a matrix of Fourier coefficients. As the change of coordinates under the requirement (23) is 1 to 1, then the matrix of the system (28) is invertible and the system always admits one solution.
Moreover, it is easy to predict what will be the error in the reconstruction of any function *f(x)* employing a specific finite number N of components; and the technique developed here could be used to ease the computation of the inverse frame operator as the limit of the system (28). Additional discussion on this method of decomposition viewed as a frame can be found in the Appendix.

### *"Indirect" method*

This method takes explicit advantage of a change of coordinates. It is based on the assumption that the reconstructed function in terms of the new basis must have the same Fourier components as the original *f(x)*; so what we will look for is the zeroing of the difference, one (Fourier) harmonic at a time.
The benefits of the method lay in the possibility of using an iterative algorithm, thus avoiding the difficulties in the inversion of the operator (the system 28).
Let us start with two real functions $S(x)$ and $R(x) \in L^2$ and, for simplicity, let us consider a function *f(x)* with a coefficient $C_0 = 0$.
If the functions $S(x)$ and $R(x)$, can be used to univocally build series that can converge to *any* function $f(x) \in L^2$, then the couple $S(x)$, $R(x)$ is a basis in that space.
The difference from $S(x)$, $R(x)$ basis, and the completely orthogonal sine, cosine basis is that the reconstruction of *f(x)* at any frequency k, as a linear combination: $A_k S(kx) + B_k R(kx)$, introduces "noise components" at frequencies multiples of k.
For this reason we are forced to use simultaneously the analysis and reconstruction phases.
As a matter of fact we could better call this procedure the "de-construction" of a function.
Let us start at the fundamental frequency. Remembering that the Fourier decomposition is:

$$f(x) = \sum_{k=1}^{\infty} \left( a_k \sin(kx) + b_k \cos(kx) \right) \qquad (31)$$

We want the Fourier components of the *f(x)* and of the reconstructed function to be equal at the fundamental frequency:
$a_1 \sin(x) + b_1 \cos(x) = [A_1 S(x) + B_1 R(x)]_1$ (32)
It can be written as:

$$a_1 \sin(x) + b_1 \cos(x) = A_1 s_1 \cos(x) + A_1 s_1^{'} \sin(x) + B_1 r_1 \cos(x) + B_1 r_1^{'} \sin(x) +$$
$$+ A_1 \sum_{q=2}^{\infty} \left[ s_q \cos(qx) + s_q^{'} \sin(qx) \right] + B_1 \sum_{q=2}^{\infty} \left[ r_q \cos(qx) + r_q^{'} \sin(qx) \right] \qquad (33)$$

The second part of this equation is the "basis noise" at higher frequencies.
As sine and cosine are orthogonal, we can separate those terms; so at frequency 1, we have:

$b_1 \cos(x) = (A_1 s_1 + B_1 r_1) \cos(x)$ (34)
$a_1 \sin(x) = (A_1 s'_1 + B_1 r'_1) \sin(x)$ (35)

We get a system of two equations with two unknowns, $A_1$ and $B_1$:

$$b_1 = A_1 s_1 + B_1 r_1 \tag{36}$$
$$a_1 = A_1 s'_1 + B_1 r'_1$$

That can be solved if the independence condition (23) is satisfied.
At the next steps of the analysis we need to take into account the "basis noise" introduced by the reconstruction at lower frequencies.
Then the input to the computation algorithm shall be the error function $f_{e1}(x)$, i.e. the function itself minus the previous approximation (note that this is the only difference with the Fourier method that uses the same $f(x)$ instead):

$$f_{e1}(x) = f(x) - A_1 S(x) - B_1 R(x) \tag{37}$$

At the next step of the analysis, for frequency k=2, the system becomes:
$$b_2 - A_1 s_2 - B_1 r_2 = A_2 s_1 + B_2 r_1$$
$$a_2 - A_1 s'_2 - B_1 r'_2 = A_2 s'_1 + B_2 r'_1 \tag{38}$$

The only unknowns here are $A_2$ and $B_2$, and the system still has one solution if the condition (23) holds.
Finally at any frequency **n** the system to solve is:

$$+ b_n - \sum_{k|n}^{n-1} A_k s_{(n/k)} - \sum_{k|n}^{n-1} B_k r_{(n/k)} = A_n s_1 + B_n r_1$$

$$+ a_n - \sum_{k|n}^{n-1} A_k s'_{(n/k)} - \sum_{k|n}^{n-1} B_k r'_{(n/k)} = A_n s'_1 + B_n r'_1 \tag{39}$$

Where: $k|n$ indicates that only the terms for which k is a factor of n are part of the equations (all the integers k that are submultiples of n, from 1, to n-1).
Then, starting from the fundamental frequency, one can calculate iteratively the coefficients $A_n$ and $B_n$.
In case a plain old sine-cosine base is employed, only the $s_1$ and $r'_1$ coefficients differ from zero, and this algorithm becomes the standard Fourier analysis.
A software routine has been written and tested in terms of analysis and reconstruction with different functions $f(x)$ and different bases $[S(x), R(x)]$. It works satisfactorily, as we can get an error as small as we want, between the $f(x)$ and the reconstructed series $S_N(x)$.

## *Results*

Both methods, the direct and indirect, have been tested in their capability to analyze and reconstruct a function $f(x)$ using a couple of basis functions $S(x)$ and $R(x)$ that satisfy the above "reconstruction requisite" (23) and "convergence requisite" (24).
However, some aspects should be taken into account.
First, even if $f(x)$, $S(x)$ and $R(x)$ all have a finite Fourier spectrum, the exact reconstruction of $f(x)$ as a series of $S(x)$ and $R(x)$ can require infinite components; nevertheless, at each step in the approximation, the difference between the function $f(x)$ and the approximation decreases in amplitude and is moved to higher and higher frequencies.
If the $S(x)$ and $R(x)$ functions do not satisfy the convergence requirements (24), then the reconstruction at any frequency can still be done, but at the cost of an increasing high frequency noise (in this sense the procedure diverges). But, as explained in the Appendix, one can always use different couples of basis functions at different frequencies.

In all the following figures, the original function *f(x)* is dotted and the reconstructed function is solid; in the bottom left corner there is a compressed graph of the basis functions and the order N of the series is indicated. The "indirect" iterative algorithm has been employed, although the difference between the two methods would be unnoticeable.

In Fig. 2 a function has been analyzed / reconstructed using a square wave and a saw tooth. The approximation is surprisingly good, given the low number of components and the "esoteric" basis employed. From an engineering point of view, a standard low pass filtering could clean out the high frequency noise, leading to a low error between the reconstructed and the original function.

In Fig. 3 the same *f(x)* has been analyzed / reconstructed using a different basis; here the approximation is better than in the previous figure, due to the similarity between the basis and the function itself.

In Fig. 4 the same analysis / reconstruction has been extended to comprise more components; now the error is lower and moved to higher frequency with respect to that of Fig. 3, as one could reasonably expect from a converging series. The noise starts where the approximation stops, in this case it goes from frequency 41 up.

An useful tool in Fourier analysis is the frequency spectrum, it is governed by the Parseval's equality:

$$\text{power(s)} = C_0^2 + \frac{1}{2} \sum_k (a_k^2 + b_k^2) \tag{40}$$

But the (51) is valid only because the sine-cosine basis is completely orthogonal.
When employing generic bases we can only use the triangle inequality to find a quantity that gives us an idea of how the coefficients $A_k$ and $B_k$ change with the frequency:

$$\text{power(s)} \leq C_0^2 + \frac{1}{T} \sum_k \int_0^T [A_k S(kt) + B_k R(kt)]^2 dt \tag{41}$$

Of course one can rewrite a new, generalized, form of the Parseval's equality.
We can plot the result of the above integral according to the frequency; it gives us an idea of the contribution of each component. In Fig. 5 there is the plot of this "generalized" spectrum relative to the functions in Fig. 4.

At this point the next natural step is to introduce a form of "generalized" filtering.
For example, from the "generalized" spectrum in Fig. 5 we can cut the first 3 components. Once reconstructed, the function corresponding to the remaining coefficients (Fig. 6) is the equivalent of a band-pass generalized filtering, in the vectorial space of the basis [S(x), R(x)]. It is a band-pass filtering because it is related only to the components from frequency 4 to 40.
Of course it has nothing in common with standard (sinusoidal) filtering; for instance, the law of energy conservation can now be violated.

## *Appendix*

Only the necessary conditions for the convergence will be given here, while the broader sufficient requisites are yet an open question. To check whether two functions S(x) and R(x) form a basis we need to verify if, given a certain reconstruction algorithm, we can find unique coefficients and that the error in norm between the original function and its reconstructed counterpart tends to zero. We split the proof into two parts. First, we check the existence and the uniqueness of the coefficients, then find the conditions that make the error tend to zero.
It is easy to verify that a couple of functions:

$$S(x) = \sum_{q=1}^{\infty} \left[ s_q \cos(qx) + s'_q \sin(qx) \right] \tag{A1}$$

$$R(x) = \sum_{q=1}^{\infty} \left[ r_q \cos(qx) + r'_q \sin(qx) \right] \tag{A2}$$

can reconstruct any (Fourier) harmonic of the function *f(x)*:

$$f(x) = \sum_{i=1}^{\infty} a_i \sin(ix) + b_i \cos(ix) = \sum_{n=1}^{\infty} A_n S(nx) + B_n R(nx) \tag{A3}$$

if their first coefficients are independent:
$|s_1 r_1'| \neq |s_1' r_1|$  (A4)
As shown in the paragraph on the Indirect Method of analysis.

Once we proved that the "independence requirement" (A4) is sufficient for a couple of functions S(x) and R(x) to assemble any sine and cosine linear combination, we need to find out the "convergence condition" that lets the series converge to any $f(x) \in L^2$ .
Let us start with the first approximation of the series. Then the error function **f_{e1}(x)** will be the original function minus the approximation relative to the fundamental frequency.
What we require is that:
$\|f_{e1}(x)\| < \|f(x)\|$  (A5)

But:

$$f_{e1}(x) = f(x) - A_1 S(x) - B_1 R(x) = (a_1 - A_1 s'_1 - B_1 r'_1)\sin(x) + (b_1 - A_1 s_1 - B_1 r_1)\cos(x) + \\ + \sum_{i=2}^{\infty} a_i \sin(ix) + \sum_{i=2}^{\infty} b_i \cos(ix) - \sum_{i=2}^{\infty} (A_1 s'_i + B_1 r'_i)\sin(ix) - \sum_{i=2}^{\infty} (A_1 s_i + B_1 r_i)\cos(ix) \tag{A6}$$

The first two terms, the coefficients of sin(x) and cos(x), are zero because $A_1$ and $B_1$ are the results of the system (36).
We now take the norm: $(<f_{e1}(x), f_{e1}(x)>)^{1/2}$ , written also as: $\|f_{e1}(x)\|$:

$$\| f_{e1}(x) \| = \| \sum_{i=2}^{\infty} a_i \sin(ix) + \sum_{i=2}^{\infty} b_i \cos(ix) - \sum_{i=2}^{\infty} (A_1 s'_i + B_1 r'_i)\sin(ix) - \sum_{i=2}^{\infty} (A_1 s_i + B_1 r_i)\cos(ix) \| \tag{A7}$$

While the norm of the *f(x)* can be written:

$$\| f(x) \| = \| a_1 \sin(x) + b_1 \cos(x) + \sum_{i=2}^{\infty} a_i \sin(ix) + \sum_{i=2}^{\infty} b_i \cos(ix) \| \tag{A8}$$

But we can substitute the first two terms from the system (36):

$$\| f(x) \| = \| (A_1 s'_1 + B_1 r'_1)\sin(x) + (A_1 s_1 + B_1 r_1)\cos(x) + \sum_{i=2}^{\infty} a_i \sin(ix) + \sum_{i=2}^{\infty} b_i \cos(ix) \| \tag{A9}$$

Comparing (A7) and (A9) we find that $\|f_{e1}(x)\| > \|f(x)\|$ when, for any $A_1$, $B_1$ not both zero:

$$(A_1 s_1 + B_1 r_1)^2 + (A_1 s_1' + B_1 r_1')^2 > \sum_{i=2}^{\infty}(-A_1 s_i - B_1 r_i)^2 + \sum_{i=2}^{\infty}(-A_1 s_i' - B_1 r_i')^2 \qquad (A10)$$

We can iterate the procedure, each time using as a new function the error function $f_e(x)$ found from the previous iteration. As at any step we annihilate one of the harmonics, the norm of the error function, which is always positive, tends to zero as the number of iterations tends to infinite.

There is also a more intuitive way of saying it if we remember that the energy of a signal is the square of the norm and is proportional to the sum of the squares of its Fourier coefficients. From the (A10) it turns out that the energy of S(x) and R(x) has to be mostly at the fundamental frequency; so any time we cancel the k-th harmonic from the error function, subtracting $A_k S(kx) + B_k R(kx)$, we add energy (as "noise" at higher frequency) lower than the energy that we subtract.
To prove that [S(x), R(x)] is a basis in $L^2$ we need to prove also that the reconstruction is unique, but this comes as a direct consequence of the use of the Fourier coefficients and the uniqueness of the solutions of the systems (36) and (39).

*An in-depth look at the convergence requirements.*
The condition (A10) imposes a requisite on the entire function energy in order to ensure the convergence but it could be even more relaxed.
To study the convergence conditions of the series of S(x) and R(x), we should look at the behavior of the $A_n$ and $B_n$ coefficients as **n** tends to infinite, or when: $a_n$ and $b_n \cong 0$.
A preliminary analysis suggests that the (necessary and sufficient) convergence condition for the series could be more relaxed than the requirement (A10), even though a thorough study is still necessary.

*Some additional considerations.*
1) If the two functions S(x) and R(x) do not satisfy the convergence criteria (A10), the approximation at each frequency can still be done, but at the cost of an increasing noise at higher frequency of the reconstructed function. In this sense, the series diverges or, better, the couple S(x), R(x), is complete but not a basis. But:

2) The new decomposition paradigm defined above allows us to use different bases at different frequencies.
In fact, one could do the analysis using a certain basis [S1(x), R1(x)] for the fundamental and up to a frequency M; then, switch to a new basis [S2(x), R2(x)] and continue the analysis from the frequency M up.
This could be used to solve any convergence problem. One could, for example, use a couple of "non-converging basis" functions for the low frequency components and allow the error function to build up high frequency noise. Then, starting from a certain frequency, one could switch to a different, converging basis, for example the ever trustworthy sinusoid, to iron out all the noise. A corollary of this statement is that any couple of independent functions can be part of a basis, as long as it is encapsulated inside a "converging" basis; this is the rationale for using the oxymoron: "diverging basis" instead of the more precise: "complete, but not a basis". An equivalent, alternative approach could take advantage of analysis-synthesis algorithms that include some kind of low pass filter, such that, the more one goes up in frequency, the more the basis gets reduced to simple sinusoids.

3) Finally, the Direct Method could be used to define the convergence requirement in terms of the characteristics of the functions S(x) and R(x) themselves, without any reference to their Fourier components.

As a matter of fact the system (31) can be viewed as a special class of the mathematical structures called Frames. A set of functions $\{f_i\}_{i \in I}$ is a Frame if, for any function $f$, there are two constants A, B > 0 called the bounds of the frame, such that:

$$A \| f \|^2 \leq \sum_{i \in I} |\langle f, f_i \rangle|^2 \leq B \| f \|^2 \qquad (A12)$$

If A=B it is called a tight frame; if $\|f_i\| = \|f_j\|$ for all i, j $\in$ I (I is an index set) it is said to be uniform; if $\|f_i\| = 1$ it is called normalized frame.

It can be seen that the structures introduced here are frames, complete, normalizable, non-tight, and studied accordingly.

**Acknowledgments**

The author wishes to thank F. Giordano of Università Parthenope Naples, S. Cavaliere and I. Ortosecco of Università Federico II Naples, P. Parascandolo of INFN for their continue support and friendship. He wishes also to thank P. Corbo and G. Langwagen of Universidad ORT of Montevideo, Uruguay for the opportunity to teach there.

**Bibliography**

[1] Daubechies I. *The wavelet transform, time-frequency localization and signal analysis*. IEEE Transactions on Information Theory, 36(5):961-1005, 1990.

[2] Daubechies, I. *Ten Lectures on Wavelets.* SIAM, Philadelphia. (1992).

[3] Daubechies, I. *From the original framer to present-day time-frequency and time-scale frames*. The J. Fourier Anal. Appl. 3 1997

[4] Duffin, R. J., Schaeffer, A. C. *A class of nonharmonic Fourier series*. Trans. Amer. Math. Soc. Vol. 72 1952

[5] Haar A. *Zur Theorie der orthogonalen Funktionensysteme*. Mathematische Annalen, 69:331-371, 1910

[6] Jackson, D. (1941). *Fourier series and orthogonal polynomials.* Mathematical Assoc. of America, Washington, D.C.

[7] Koc, C. K., Chen, G., Chuy, C. K. *Analysis of computational methods for wavelet signal decomposition and reconstruction*. IEEE Transaction on Aerospace and Electronic Systems. July 1994

[8] Mallat, S. G., *A theory for multiresolution signal decomposition: the wavelet representation*. IEEE Transactions on Pattern Recognition and Machine Intelligence. Vol. 11, No. 7, pp. 674-693, July 1989.

[9] Walker, J.S. : *Fourier Series*  http://www.uwec.edu/walkerjs/media/fseries.pdf

[10] Walker, J.S. (1988). *Fourier Analysis.* Oxford Univ. Press, Oxford.

[11] Wei, Y., Chen, N.: *Square wave analysis*. Journal of Mathematical Physics Vol. 39, N. 8, August 1998.

[12] Wei,Y., *Frequency analysis based on general periodic functions*. Journal of Mathematical Physics Vol 40 N. 7, July 1999

[13] Wei, Y.: *Frequency analysis based on easily generated functions*. Applied and Computational Harmonic Analysis, 2000

[14] http://www.siam.org/journals/problems/downloadfiles/01-002.pdf

[15] Wittgenstein L., *Tractatus Logico-Philosophicus.* 1921